\documentclass[letterpaper, margin=1in]{article}
\usepackage{amsmath, amsthm, amssymb, geometry}    
\usepackage{graphicx}			
\usepackage[colorlinks=true, urlcolor=blue, citecolor=black, linkcolor=black]{hyperref}  
\usepackage{subcaption}			

\title{A Braid Box}

\author{Blake Karl Winter\textsuperscript{1} and Amanda Taylor Lipnicki\textsuperscript{2}
\vspace{10pt}\\
\textsuperscript{1}Math Division, University of Findlay, Ohio, USA; bw3073@findlay.edu\\
\textsuperscript{2}Alfred University, New York, USA; tayloral@alfred.edu} 

\date{}					

\begin{document}

\maketitle

\thispagestyle{empty}

\begin{abstract}

We give a method for constructing an interactive art piece which illustrates two different definitions of the braid groups, along with their faithful action on the free group. The box also demonstrates how all motions of points in the plane can be realized by motions in a single T-shaped subspace of the plane. This helps students and those who are not specialists in algebraic topology to understand these important topological objects.

\end{abstract}


\section*{The braid groups}
The braid groups are extremely useful for the study of topology and algebra. First considered by Artin \cite{Artin}, they were later the focus of study by Birman and others \cite{Birman}. While braids can be drawn in diagrams to describe them to non-specialists, there are other equivalent definitions which are more difficult to explain or intuitively grasp. Herein, we give a method for constructing a box which illustrates how two different definitions of the braid groups are related. First, we consider the definitions of the braid group which we will be using.

\begin{figure}[ht!]
	\centering
    \begin{minipage}{0.45\textwidth}
    \begin{flushright}
	\includegraphics[width=2.5in]{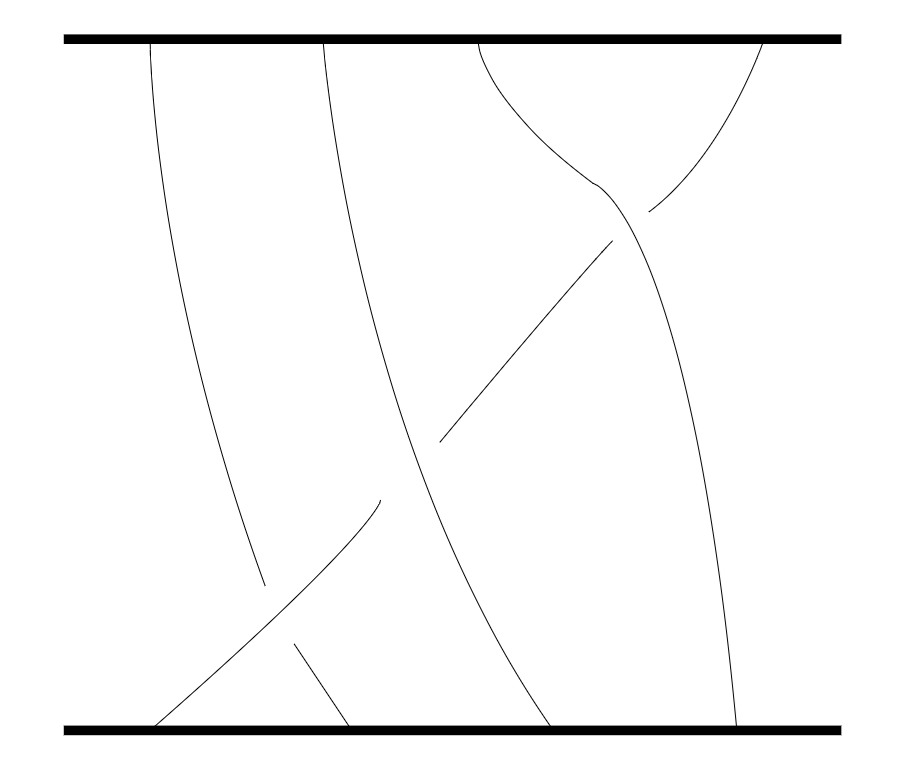}
    \caption{A 4-strand braid.}
    \label{4-strand}
    \end{flushright}
    \end{minipage}
    \begin{minipage}{0.45\textwidth}
    \includegraphics[width=2.5in]{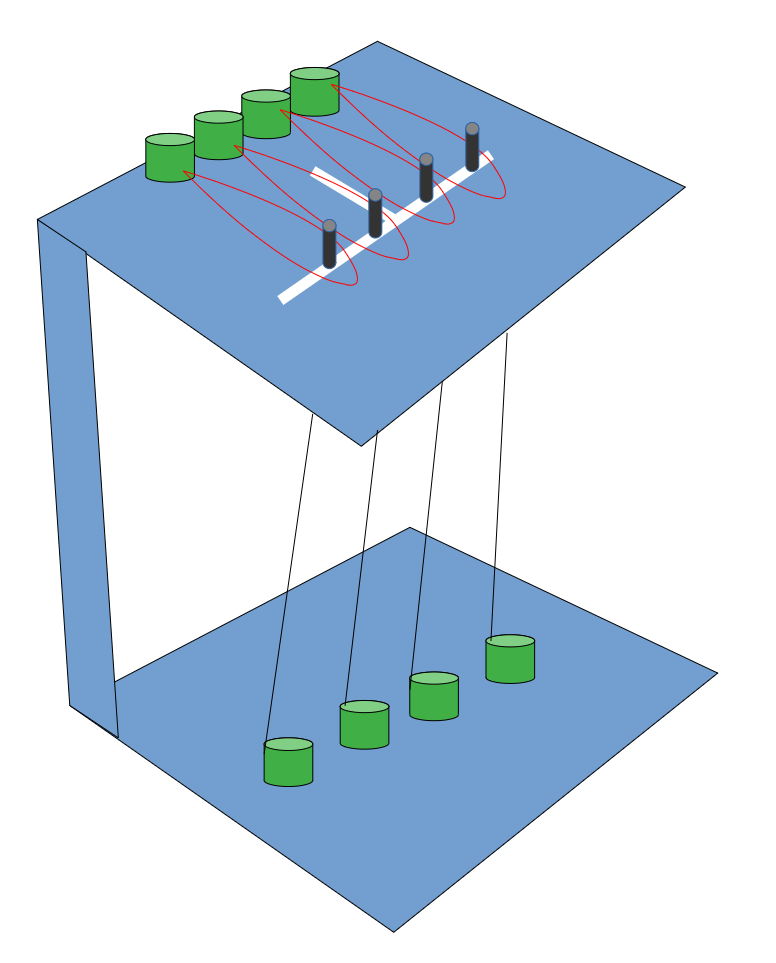}
    \caption{An illustrated schematic of the box design with 4 strands.}
    \label{box}
    \end{minipage}
\end{figure}

One definition of an $n$-braid is proper embeddings $f$ of the product $\sqcup_{i_1}^n I$ into $\mathbb{R}^2 \times I$, where the tangent vector of each copy of the interval always has a positive component along projection to the third axis. Such an embedding, in general position, can be given by a \emph{braid diagram} like the one illustrated in Figure \ref{4-strand}.

Equivalently, we may take proper embeddings $f$ of the product $\sqcup_{i_1}^n I$ into $\mathbb{R}^2 \times I$, where the embedding meets level sets of the $I$ component of $\mathbb{R}^2 \times I$ transversely. This means that the projection to the $I$ component defines a Morse function on the embedding with no births or deaths.

Two such braids are considered equivalent if they are related by an (ambient, PL, or smooth) isotopy that leaves the boundary fixed. The group multiplication is given by stacking one braid on top of another. Under this definition, the braid group on $n$ strands, $B_n$, has the following presentation:

$B_n = \langle b_1, b_2, ..., b_n | b_i b_{i+1} b_i = b_{i+1} b_i b_{i+1}, b_i b_j = b_j b_i, |i-j| \neq 1 \rangle $.

Another definition may be given in terms of a motion group. Let $C_n$ be the configuration space of $n$ points in the plane. Then a motion is a loop in $C_n$. Two braids are equivalent if the two loops are homotopic. 

\begin{figure}[ht!]
    \centering
    \includegraphics[width=0.5\linewidth]{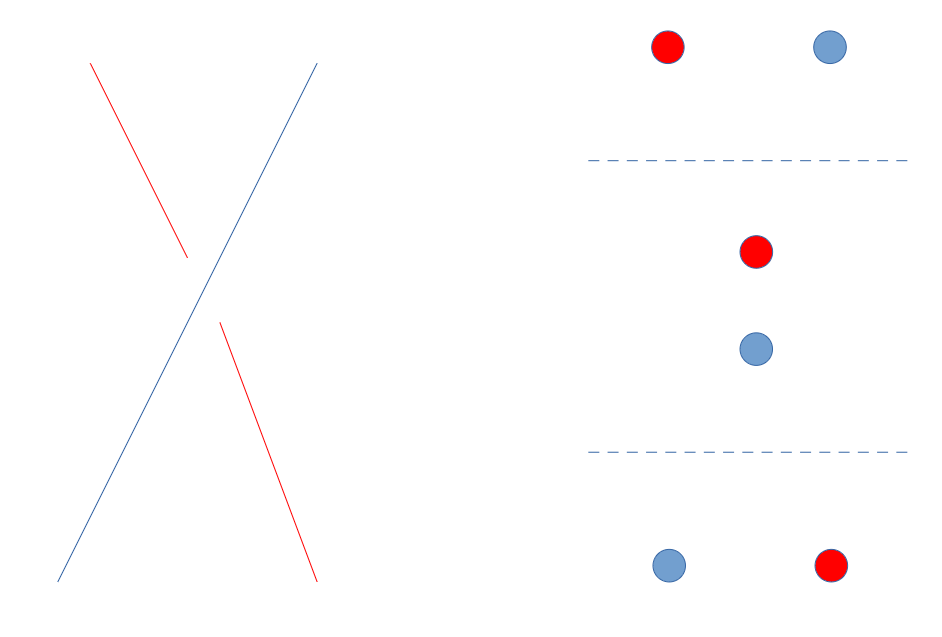}
    \caption{The correspondence between a braid and a motion of points in the plane.}
    \label{Motion}
\end{figure}

The relationship between these two definitions is topologically straightforward---any function $f:I \rightarrow C_n$ defines an embedding of a disjoint union of intervals into $\mathbb{R}^2 \times I$ by sending the point $t$ of the $i$-th component to the position of the $i$-th point in the configuration space at $f(t)$. In other words we view the braid as foliated along the $t$ direction, and each constant-$t$ slice corresponds to a configuration in $C_n$. This may be reversed as well, so that given a braid, we foliate it along the $t$ direction, and then define a loop in $C_n$ such that the configuration at any given $t$ corresponds to the configuration of the slice of the braid at that $t$. This constructs an inverse map. Two braids are equivalent if and only if the loops in $C_n$ are homotopic, as observed by Artin \cite{Artin}. Therefore, the $n$-braid group is isomorphic to the fundamental group of $C_n$. An example of a simple braid and the corresponding motion of points is given in Fig. \ref{Motion}.

The relationship between these definitions can be clarified by thinking of foliating $\mathbb{R}^2 \times I$ along the second component. Then each slice in the foliation involves a configuration of $n$ points. The braid then traces out a motion. Or, conversely, we may see a motion as dragging a braid behind it and braiding it as it goes. In our box, this relationship is explicitly demonstrated, as moving the dowels (representing the points) pulls strands beneath them into the form of the braid corresponding to the motion.

\begin{figure}[ht!]
    \centering
    \begin{minipage}{0.45\textwidth}
    \begin{flushright}
    \includegraphics[width=0.8\linewidth]{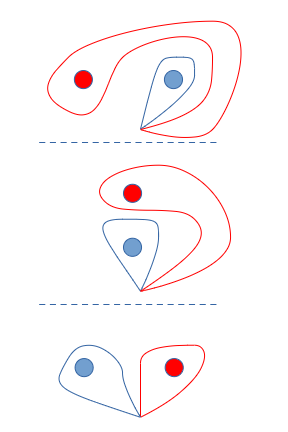}
    \caption{A diffeomorphism that moves the points in the plane also moves the generators of the fundamental group of their compliment in the plane, which is a free group. This induces a faithful action on the free group $F_n$ by the $n$-strand braid group $B_n$.}
    \label{Action}
    \end{flushright}
    \end{minipage}
    \begin{minipage}{0.45\textwidth}
    \includegraphics[width=\linewidth]{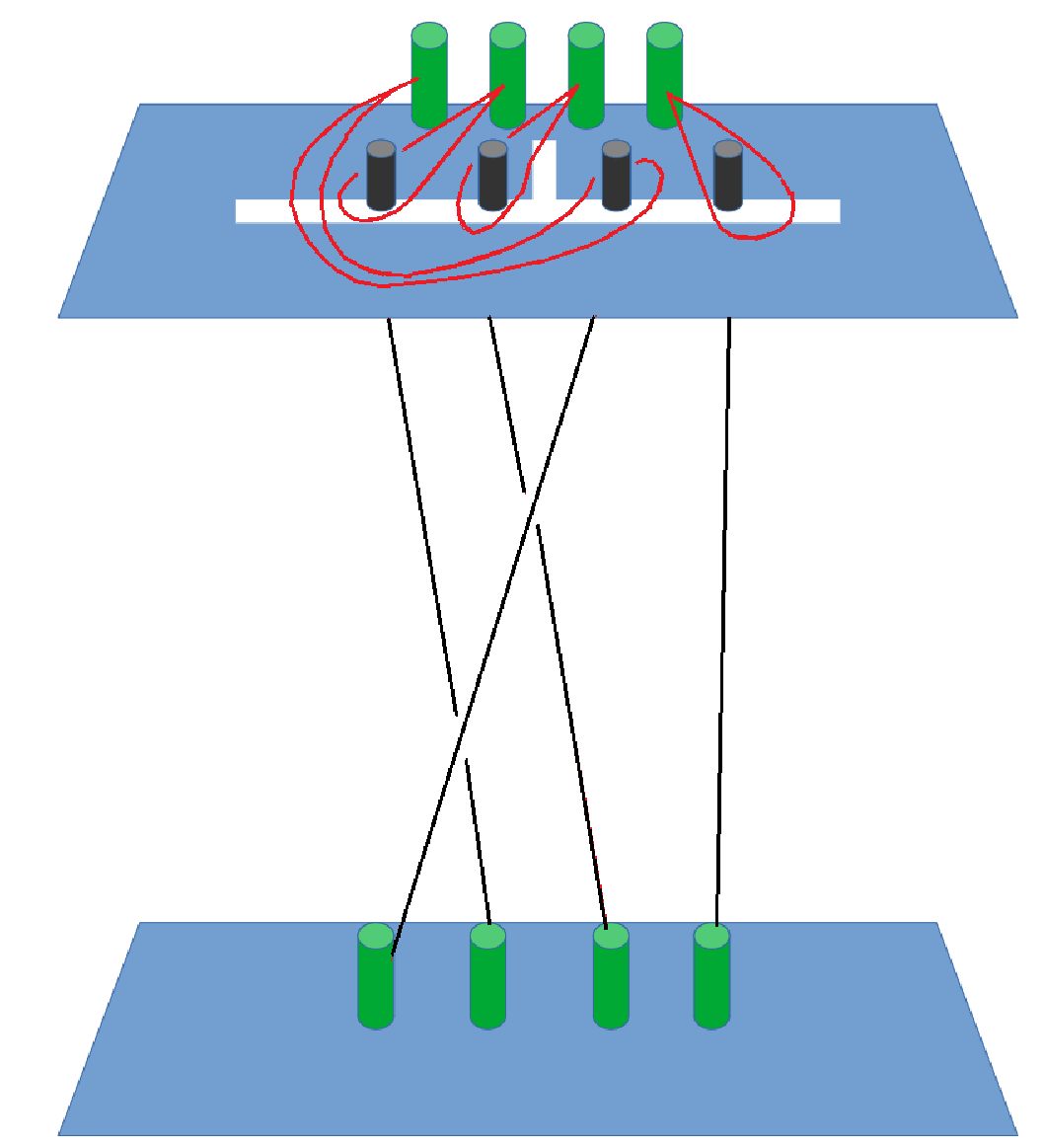}
    \caption{The box after moving the third dowel behind the first and second dowels. The resulting braid is shown, along with the action on the free group.}
    \label{boxuse}
    \end{minipage}
\end{figure}

There is another mathematical property of braids which is of particular interest; $B_n$ acts on the free group on $n$ generators. This action can be seen by considering a motion of $n$ points in the plane. The motion extends to an ambient isotopy, and at the end of the motion, this determines a diffeomorphism from the $n$-punctured plane to itself. The fundamental group of the $n$-punctured plane is a free group on $n$ elements. The generators are embedded loops from some chosen basepoint, each one wrapping around a single point in the plane once. Thus, a diffeomorphism from the $n$-punctured plane to itself induces an automorphism of the fundamental group, by rearranging the loops that generate it. Indeed, the automorphism will be an inner automorphism. Our box also demonstrates this. The action itself is illustrated in Fig. \ref{Action}. In our box, loops passing around the dowels represent the generators of the fundamental group, and as the dowels are moved, the loops move to illustrate the action.

Finally, it is a curious topological fact that all loops in the configuration space $C_n$ are homotopic to loops which only pass through configurations in which the $n$ points lie in $[-1, 1]\times \{ 0 \} \cup \{ 0 \} \times [0, 1]$. In general the points move back and forth along the first axis. If one point passes around another, we move  it up the vertical axis, pass the other point in front of it, and then bring it back. This fact allows us to build a box which is physically realizable.

Such a description of a braid (or link) is known as a \emph{3-page book} description \cite{Bernhart}. This method of describing links and braids is not particularly common, but it is known that all of  them can be placed into such a 3-page book. This box illustrates this fact for braids. In particular, the motion of the dowels is constrained to lie within a T shape, but (up to topological equivalence) any motion/braid may still be generated (provided that there is enough string).

\section*{The box}

The box consists of a lower platform with $n$ spools. The spools have strings which attach to movable dowels. There is an upper platform with a T-shape cut into it. The dowels have disks in their middles so that they can be inserted into this groove but not fall through.

On the bottom of each dowel is a hook to which a string can be attached. As the dowels are slid around on the top platform, they perform a motion of $n$ points in the plane. Underneath, the strands are pulled into the corresponding braid. The braid is kept reasonably taut by the spools, making the braid easy to see. After creating a braid, the viewer can  return the setup to its original configuration by unhooking the strands, spooling them back up, and then reattaching them to the dowels, which is useful in cases where the braid is sufficiently complicated that untangling it becomes difficult. However, showing how the inverse motion may be read off the braid is also possible with the help of the box.

On the top of the upper platform, there are an additional $n$ spools with loops. Each loop may be placed over the top of a dowel. These loops represent the generators of the fundamental group of the plane with $n$ points removed. As the dowels are moved to form a braid, the loops are pulled with them, showing the automorphism which the braid generates. See Fig. \ref{box} for an illustration of a braid box with 4 strands. A different angle, with a nontrivial braid and the corresponding action, is shown in Fig. \ref{boxuse}.

\section*{Variations}

One possible variant would be to have the loops, corresponding to the generators of the fundamental group, placed around the strands between the two platforms. Then the loops could be pushed up and down the braid itself, showing how the braid moves the generators.

It is also possible to illustrate the action on the fundamental quandle if the loops on the top side are replaced by strands that hook to the dowels. The quandle is a very useful algebraic object for the study of knots and links, introduced by Joyce \cite{Joyce}. The generators of a quandle are topologically represented by curves from a basepoint to the boundary of a tubular neighborhood of the link (or braid). Such a curve may be physically represented by a strand that attaches to the dowel. The quandle operation is a self-distributive operation. Given two paths $a, b$ from the basepoint  to the boundary of a tubular neighborhood of the link, the path corresponding to $a*b$ is created by following $b$, then a meridional loop around the link at the endpoint of $b$, then following $b$ in reverse, and finally following $a$. Paths are considered up to homotopy that keeps their initial point fixed at the basepoint and their terminal point on the boundary of the tubular neighborhood of the link.

In future, we may consider constructing artistic illustrations of other motion groups, such as the loop braid group.

\section*{Conclusion}

The braid groups are an important topic in mathematics with simple diagrams appealing to non-specialists. However, some of the more technical aspects of their definition and their behavior are difficult to understand. This box provides an interactive way to illustrate these relationships while also providing a beautiful audience-engagement driven art piece.

    
{\setlength{\baselineskip}{13pt} 
\raggedright				

} 
   
\end{document}